\documentclass{amsart}

\usepackage{amssymb}
\usepackage{graphicx}
\usepackage{import}
\usepackage{tikz}
\usepackage{pgfplots}
\usepackage{mygraphs}
\usepackage{cite}

\theoremstyle{definition}

\numberwithin{equation}{section}

\newcommand{\setbar}{\,|\,} %bar for sets

\title{Loop-Erased Random Surfaces}

\author{Kyle Parsons}
\address{Department of Mathematics, The Ohio State University, 
Columbus, OH 43210}
\email{parsons.299@math.osu.edu}

\begin{document}

\begin{abstract}
Loop-erased random walk and it's scaling limit, Schramm--Loewner evolution, have found numerous applications in mathematics and physics.  We present a 2 dimensional analogue of LERW, the loop erased random surface.  We do this by defining a 2 dimensional spanning tree and declaring that LERS should have the same relation to these 2 trees as LERW has to ordinary spanning trees.  Furthermore we present numerical evidence that the growth rate for LERS on a $\delta$ fine grid as $\delta \to 0$ is $2.5269 \pm 0.0017$ and we hypothesize that it has an exact value of 48/19.  This suggests the possibility of a fractal limiting object for LERS analogous to SLE for LERW.
\end{abstract}

\maketitle

%%%%%%%%%%%%%%%%%%%%%%%%%%%%%%%%%%%%%%%%%%%%%%%%%%%%%%%%%%%%%%%%%%%%%%
\section{Introduction}
%%%%%%%%%%%%%%%%%%%%%%%%%%%%%%%%%%%%%%%%%%%%%%%%%%%%%%%%%%%%%%%%%%%%%%
\subsection{Preliminaries}
Loop erased random walks (LERW) are a model of self-avoiding random walks introduced by Lawler in 1980 \cite{Lawler:1980}.  In that paper he proved, however, that the loop erasure of a simple random walk on the integer lattice does not result in the uniform measure on self avoiding walks.  Informally speaking, a loop erased random walk is generated by chronologically removing loops from a simple random walk. Since their introduction, loop erased random walks have been studied extensively \cite{Barlow:2010}\cite{Kenyon:2000}\cite{Schramm:2000}.  They have a close connection with uniform spanning trees that is profitably exploited by Wilson's algorithm \cite{Wilson:1996}.  Using Wilson's algorithm we can generate a uniform spanning tree on a graph, by repeatedly sampling loop erased random walks.  Conversely, the path between two points in the uniform spanning tree has the same distribution as the loop erased random walk between the points.  

On the two dimensional integer lattice a LERW starting from the origin and stopped when it first leaves a ball of radius $r$ centered on the origin will have around $r^{5/ 4}$ steps, and the scaling limit of LERW in two dimensions is Schramm--Loewner Evolution with parameter 2 ($\text{SLE}_2$), which is known to be conformally invariant and have fractal dimension $5/ 4$\cite{Kenyon:2000}\cite{Schramm:2000}\cite{Schramm:2005}\cite{Lawler:2004}.  Schramm-Loewner Evolution has found a myriad of uses in physics.  For an overview see Cardy \cite{Cardy:2005}. 

\begin{figure}
\includegraphics[width=\textwidth]{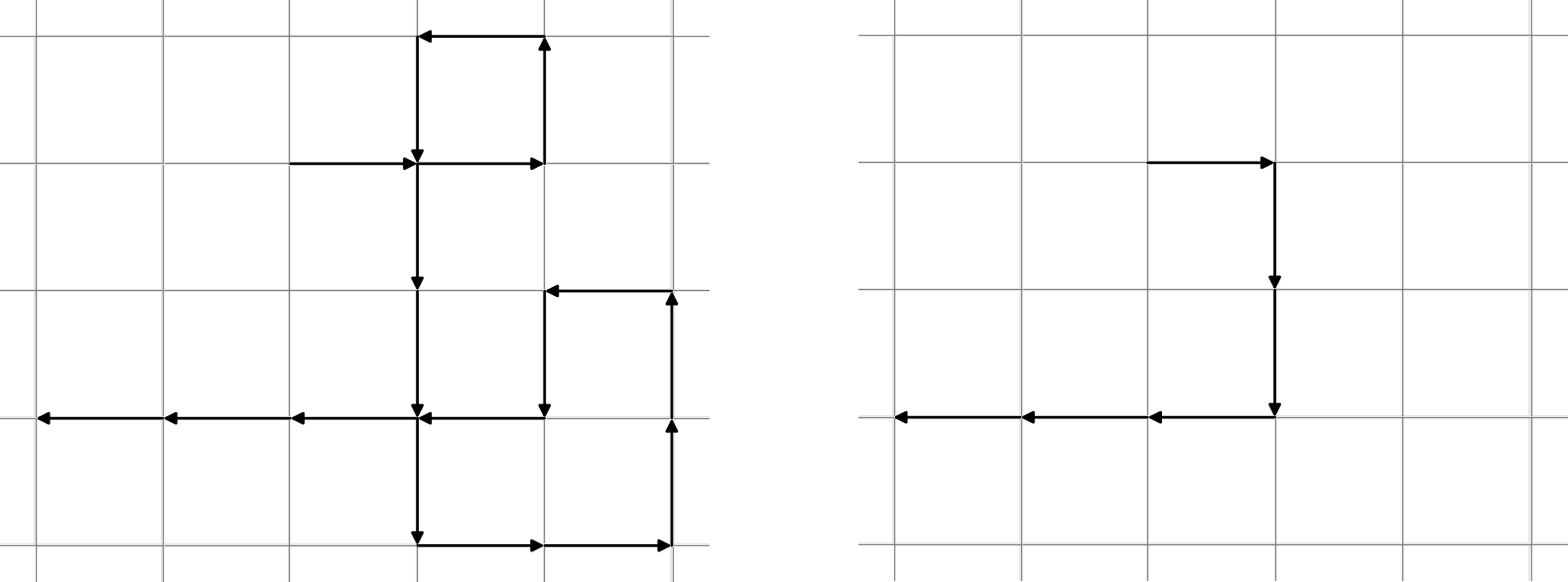}
\caption{A simple random walk on the left and its loop erasure on the right}
\end{figure}

\begin{figure}
$\begin{array}{cc}
\includegraphics[width=0.5\textwidth]{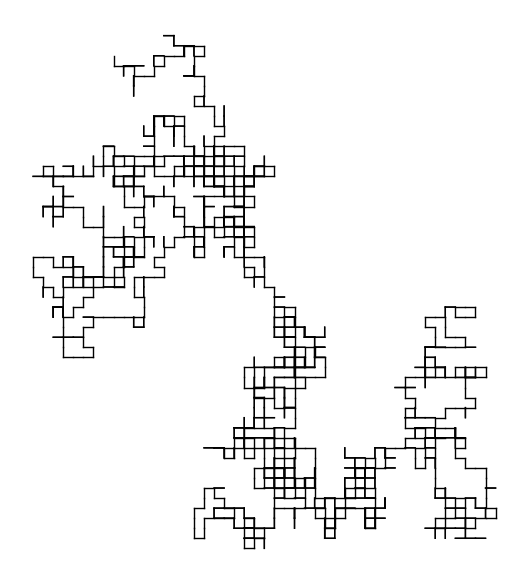} & \includegraphics[width=0.5\textwidth]{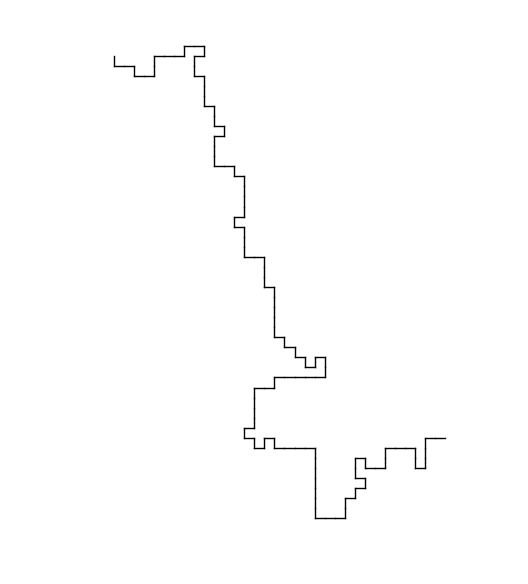} \\
\end{array}$
\caption{A simple random walk run until the first exit from a ball of radius 100 about the origin on the left and its loop erasure on the right}
\end{figure}

Several different generalizations of spanning trees to higher dimensions have appeared in the literature.  One of the first was by Kalai in 1983 \cite{Kalai:1983} when he extended Cayley's formula to enumerate certain tree-like subcomplexes of the high-dimensional simplex.  Another generalization was given recently by Hiraoka and Shirai in 2015 \cite{Hiraoka:2015} with their $k$-spanning acycles.  They showed a limiting behavior for expected weight of a minimal spanning acycle as an extension of Frieze's theorem on the minimal weight of spanning trees on a random graph.  In 2009, Lyons introduced $k$-bases and complements of $k$-cobases \cite{Lyons:2009}.  In addition he introduced a nontrivial probability measure on these tree analogues.

We will work with tree-like subcomplexes of a finite lattice of square cells.  We will call a subcomplex of our lattice a \textit{$2$-tree} if it is maximal while having no integer homology of any dimension.  Note that in this paper, $2$-trees will conincide with $2$-spanning acycles, $2$-bases, and complements of $2$-cobases.  Also, the probability measure introduced by Lyons on these objects will just be the uniform measure here.  We call the $2$-tree chosen via Lyons' probability measure the \textit{determinantal $2$-tree}.  Mirroring the pairing that Wilson's algorithm gives us between spanning trees and LERW, we define the \textit{loop-erased random surface} with boundary loop $a$, to be the unique surface bounded by $a$ inside the determinantal $2$-tree.   \cite{Aizenman:1983}\cite{Grimmett:2010}.

We will work specifically with 2-trees on a finite 3-dimensional lattice of plaquettes.  In this case the determinantal measure is just the uniform measure.  We use this combined with a duality found in \cite{Lyons:2009} to generate these bases quickly using the Aldous--Broder algorithm \cite{Aldous:1990}\cite{Broder:1989}.  Concurrently, we also calculate the loop-erased random surface with boundary the equatorial loop.  In this way, generating a loop-erased random surface takes $\Theta(n^3\log(n))$ time \cite{Johan:2000} where $n$ is the length of each side of the lattice.  This avoids linear programming or any other standard techniques for finding bounded surfaces.  See for example \cite{Dunfield:2011}.

\begin{figure}
\includegraphics[width=\textwidth]{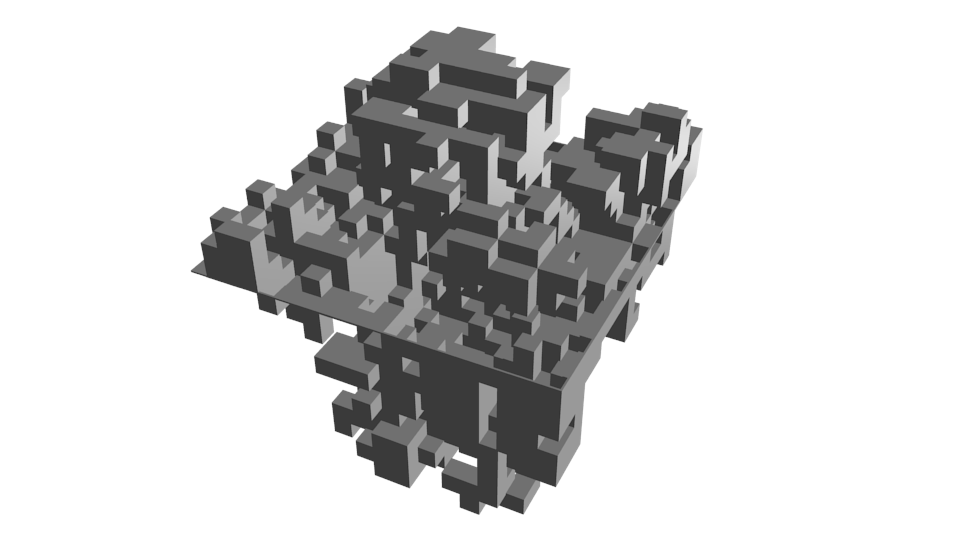}
\caption{A loop erased random surface with a 20$\times$20 loop as its boundary}
\end{figure}

\subsection{The Analogy to Graphs}
Note that throughout this paper we will use the term \textit{acyclic} to refer to a complex with trivial $\mathbb{Z}_2$ homology, and correspondingly, chain groups and homology groups will have $\mathbb{Z}_2$ coefficients unless otherwise noted.  Given a connected graph $G$, a spanning tree $H$ of $G$ is subgraph of G containing all the vertices of $G$ which is maximal with no cycles (no first homology).  Equivalently we can say that spanning trees are minimal connected subgraphs or that they are connected acyclic subgraphs.  In a similar way given a 2-dimensional complex $K$ we define a \textit{2-tree} $H$ of $K$ to be a maximal acyclic subcomplex of $K$.  

One useful property of trees is that given any two vertices of a tree there is a unique path whose ends are those two points.  We know that this path is unique because if there were two such paths, their union would form a loop.  But, by definition, trees do not have cycles.  This is especially useful since if the tree is the uniform spanning tree on the graph, then the path is the loop erased random walk between the points \cite{Wilson:1996}.  If a 2-complex $K$ is homologically simple enough, its 2-trees will span the 1-skeleton of $K$.  Then given any 1-cycle $a$ in $K$ we can define the loop-erased sufrace with a as its boundary to be the realization of a 2-chain $b$ in a 2-tree of $K$ with $\partial b = a$. For a given 2-tree this $b$ is unique since if there were two such 2-chains, their sum would be a 2-cycle.  This cannot happen because 2-trees have trivial homology. Furthermore if we choose the 2-tree from a uniform measure then we have a distribution on 2-chains bounded by $a$.  We will call this distribution the loop-erased random surface conditioned to have boundary $a$.

\subsection{Expectations and Results}
If loop erased random surfaces behave similarly to loop erased random walks, then in particular, we expect (and would like) for them to have scaling limits in the correct circumstances.  Loop erased random walks on the integer lattice have scaling limits in all dimensions.  In 2 dimensions the scaling limit is $\text{SLE}_2$.  In dimension 4 the scaling limit is Brownian motion with a logarithmic correction and in 5 or more dimensions it is simply Brownian motion with no correction needed.  For dimension 3 the scaling limit of LERW is know to not be Brownian motion and is conjectured to have dimension $13/8$ \cite{Anton:2001}.  In addition to a scaling limit for LERS giving an analogue for SLE in a low number of dimensions, it could also give us an analogue of Brownian motion in higher dimensions.  

The curve $\text{SLE}$ is conformally invariant. However, the scaling limit of LERW in 3 dimensions is known not to be (since a random walk in 3 dimensions is transient), but it is invariant under dilations and rotations \cite{Kozma:2007}.  We will see that there is a kind of duality between LERW and LERS in three dimensions.  This duality in particular leads us to expect that LERS in 3 dimensions has a scaling limit and that limit is invariant under dilations and rotations.  

Given a loop erased random surface on a grid with a specified boundary loop, as the mesh size $\delta$ of the grid goes to zero we expect that the area of our bounding surface will tend to infinity roughly at the rate $\delta^{-c}$ where the growth exponent $c$ here is of great interest.  In fact, based on numerical experiments presented here, it appears that this growth exponent exists and is approximately equal to 2.5269.

%%%%%%%%%%%%%%%%%%%%%%%%%%%%%%%%%%%%%%%%%%%%%%%%%%%%%%%%%%%%%%%%%%%%%%
\section{Notation and Definitions}
%%%%%%%%%%%%%%%%%%%%%%%%%%%%%%%%%%%%%%%%%%%%%%%%%%%%%%%%%%%%%%%%%%%%%%
Given a oriented cell complex $X$ write $\Xi_k X$ for its collection of $k$-cells and $\Xi X$ for all the cells of $X$.  Let $C_k(X,\mathbb{C})$ denote the $k$-chains of X, $Z_k(X, \mathbb{C})$ the $k$-cycles, and $B_k(X, \mathbb{C})$ the $k$-boundaries all with $\mathbb{C}$ coefficients and defined with respect to the standard cellular boundary map, $\partial_k$.  Correspondingly, $C^k(X, \mathbb{C})$, $Z^k(X, \mathbb{C})$, and $B^k(X, \mathbb{C})$ denote the cochains, cocycles, and coboundaries of the the coboundary map $\delta_k$.   Given $T\subseteq \Xi_k X$ let $X_T$ denote the subcomplex $T\cup\bigcup_{i=0}^{k-1} \Xi_k$.  By Lyons\cite{Lyons:2009}, if a subset $T\subseteq\Xi_k X$ is maximal with $Z_k(X_T)=0$ then we call $T$ a \textit{$k$-base} of X.  It is important to note that the 1-bases of connected complex are spanning trees on its 1-skeleton.  If it is maximal with $\delta_k$ injective then we call it a \textit{$k$-cobase}.  Lyons also defines probability measures $P_k$ and $P^k$ on the sets of $k$-bases  and complements of $k$-cobases respectively.  Importantly if $H_k(X)=0$ then these probability measures agree.  Also, for a connected complex $X$, $P_1$ is always the uniform measure on spanning trees of the 1-skeleton of $X$.  $P_k$ weights bases by the squared size of the torsion in dimension $k-1$.  $P^k$ weights complements of $k$-cobases by slightly different method.  In this paper, these measures will always be uniform.

We will work with a slightly different generalization of spanning trees.  Given a finte $d$-dimensional complex $X$ (simplicial, cubical, etc...) call maximal acyclic subcomplex of $X$ of dimension $k$ a \textit{$k$-tree}.  Furthermore, for a $k$-tree $T$ of $X$, if the $m$-skeleton of $T$ matches the $m$-skeleton of $X$ then we say that $T$ \textit{$m$-spans} $X$.  Given a $k$-tree $T$ of $X$ that $(k-1)$-spans $X$ and $(k-1)$-cycle $a$ in $X$, we know that $a$ is also a cycle in $T$.  Furthermore, there is a unique $k$-chain $b$ of $T$ with $\partial b=a$.  The chain $b$ exists because $T$ has no homology in dimension $k-1$ and is unique because if another such chain $b'$ existed, the sum $b+b'$ would be a $k$-cycle but $T$ is acyclic and has dimension $k$.  We call the support of the chain $b$ the \textit{loop-erased chain} with boundary $a$.  If, furthermore, the $k$-tree is chosen uniformly from the set of all $k$-trees on $X$, then we call the support of $b$, the \textit{loop-erased random chain} with boundary $a$.  In this case we call the uniform distribution on $k$-trees the \textit{uniform spanning $k$-tree} on $X$.  Because we are working with $\mathbb{Z}_2$ coefficients we can directly think of the loop-erased random chain as a subcomplex.

Let $\mathcal{Q}_n$ be the 2-dimensional cubical complex with vertex set $V(\mathcal{Q}_n) = \{(a,b,c) \in \mathbb{Z}^3:0 \le a,b,c \le n+1\}$ and whose face set consists of all unit squares supported on the vertex set.  Consider $\mathcal{Q}_n$ to be a the 2-skeleton of a regular CW structure $X_n$ on $S^3$ with $n^3 + 1$ 3-cells ($n^3$ cube shaped cells and one outer cell).  Let $X_n^*$ be the dual CW structure to $X_n$ and let $\mathcal{G}_n$ be the 1-skeleton of $X_n^*$. (That is, $\mathcal{G}_n$ is a graph with a vertex at the center of each cube of $X_n$ and one more vertex ``at infinity". There is an edge connecting two vertices of $\mathcal{G}_n$ for each face of $X_n$ that separates the vertices.)  Let $D:\Xi X_n \to \Xi X_n^*$ be the dualizing map that pairs a face $\sigma$ of $X$ of dimension $k$ with the face $D[\sigma]$ of $X_n^*$ of dimension $3-k$ which intersects $\sigma$ transversely.  We consider this map also as a map between sets of faces of $X_n$ and $X_n^*$.  We seek to uniformly sample 2-trees of $\mathcal{Q}_n$.  Given an acyclic subcomplex $T$ of $\mathcal{Q}_n$, if it does not 1-span $\mathcal{Q}_n$ it may be extended by first including every edge of $\mathcal{Q}_n$ missing and then adding plaquettes until all 1-cycles are once again boundaries.  Thus we see that all 2-trees of $\mathcal{Q}_n$ are 1-spanning and so these 2-trees are precisely the same as Lyons' 2-bases.  Given $K^*$ the complement of a 2-base of $X_n^*$ there is a coupling between $P^1$ of $X_n^*$ and $P_2$ of $X_n$ that pairs $K^*\subseteq X_n^*$ with $X_n\backslash D[K^*]$ \cite{Lyons:2009}.  Now $H_1(X_n^*)=0$ so $P^1$ is the same as $P_1$ on $X_n^*$, that is it is the uniform measure on spanning trees of $\mathcal{G}_n$.  Essentially, we uniformly generate a spanning tree of $\mathcal{G}_n$ and then remove every face of $\mathcal{Q}_n$ that is crossed by an edge of that spanning tree.  This leaves us with the desired 2-tree.  This has the benefit of there being several efficient algorithms for generating uniform spanning trees.

An important tool that we will use for sampling 2-trees will be the Aldous--Broder algorithm\cite{Aldous:1990}\cite{Broder:1989}.  The algorithm samples uniformly from all spanning trees of a finite graph.  Note that the graph need not be simple.  We have reproduced it here for reference.

\begin{itemize}
\item Choose a vertex $v_0$ by any method.
\item Let $T_0$ be the tree consisting only of $v_0$
\item Perform a random walk $v_i$ starting at $v_0$ and keeping track of the edges traversed.
\item Whenever the random walk is at a vertex $v_k$ not in the tree-so-far, $T_i$, let $T_{i+1}$ be $T_i$ along with $v_k$ and the edge $e_k$ traversed to reach $v_k$
\item Once all vertices have been visited by the random walk, return the current tree-so-far, $T_i$
\end{itemize}

%%%%%%%%%%%%%%%%%%%%%%%%%%%%%%%%%%%%%%%%%%%%%%%%%%%%%%%%%%%%%%%%%%%%%%
\section{Experiments}
%%%%%%%%%%%%%%%%%%%%%%%%%%%%%%%%%%%%%%%%%%%%%%%%%%%%%%%%%%%%%%%%%%%%%%
\subsection{Filling Area of an Equatorial Loop}

We seek to better understand the uniform spanning 2-tree on $\mathcal{Q}_n$ and to that end we studied the size of the loop-erased random surface with boundary the equatorial loop .  By size, we simply mean the number of faces in the support of the loop-erased random surface. The equatorial loop on $\mathcal{Q}_n$ is the loop with vertex set $\{(x, y, \lfloor n/2 \rfloor) \setbar (0\le x\le n \text{ and } y\in\{0,n\}) \text{ or } (0\le y\le n \text{ and } x\in\{0,n\})\}$ (that is, the loop that is the ``equator" of the cube). The loop-erased random surface with boundary the equatorial loop is analogous to the loop-erased random walk between two given points.  To that end we expect that as $n$ increases (or equivalently as our mesh gets finer), that the number of 2-cells in the LERS to grow like $n^c$ for some constant $2\le c\le 3$ known as the growth exponent.  For analogy, a loop-erased random walk on $\mathbb{Z}^2$ from the origin to a circle of radius $n$ centered on the origin will have about $n^{5/ 4}$ steps.  That is, LERW has a growth exponent of $5/ 4$.

Let $M_n$ denote the size of the loop-erased random surface with boundary the equatorial loop in $\mathcal{Q}_n$, and let $\widehat{M}_n$ denote our empirical estimate of $M_n$.  In order to estimate what the growth factor of LERS might be, we will sample from $M_n$ repeatedly and for a range of $n$, then we will find the line of best fit for this data.  The slope of that line will be our estimate of the growth exponent.  By doing this we can eliminate much of the error due to small grid size.  We will call the growth exponent of LERS $c$ (assuming it exists) and our estimate $\tilde{c}$.

In order to sample the loop-erased random surface on $\mathcal{Q}-n$ we use the following algorithm.
\begin{itemize}
\item Let $S_0$ consisting of the collection of $n^2$ squares lying in the plane spanned by the equatorial loop be an initial bounded surface for the equatorial loop.
\item Begin the Aldous--Broder algorithm on $\mathcal{G}_n$ starting at $\infty$.
\item For each edge $e_k$ added to your partial spanning tree $T_k$ by the Aldous--Broder algorithm, let $e_k^-$ be the vertex already in the $T_k$ and let $e_k^+$ be the new vertex.
\item If $e_k$ crosses the current partial bounded surface $S_i$, update $S_i$ by letting $S_{i+1} = S_i + \delta D[e_k^+]$ where we are considering our coefficient ring to be $\mathbb{Z}/2$.
\item When the Aldous--Broder algorithm finishes, return the current bounded surface $S_i$ 
\end{itemize}
At every step in this process it is easy to check that $S_i$ lies inside $\mathcal{Q}_n\\D[T_k]$ and has boundary the equatorial loop.  Furthermore, the resulting surface is independent of $S_0$ simply because it is the unique bounded chain in the resulting 2-tree.  If there were two such different chains, there sum would be a nontrivial 2-cycle which cannot exist since by definition 2-trees do not have nontrivial 2-cycles.

\begin{figure}
$\begin{array}{cc}
\includegraphics[width = 0.5\textwidth]{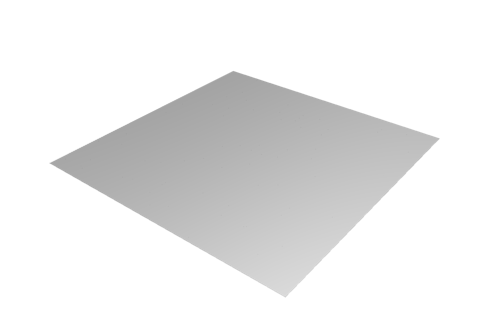}   & \includegraphics[width = 0.5\textwidth]{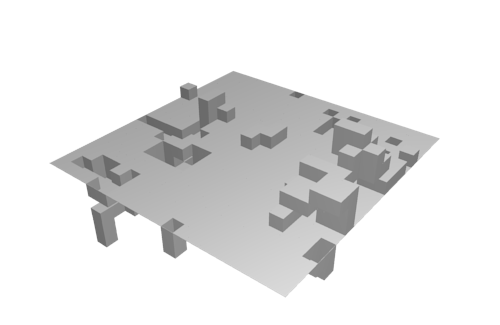} \\
\includegraphics[width = 0.5\textwidth]{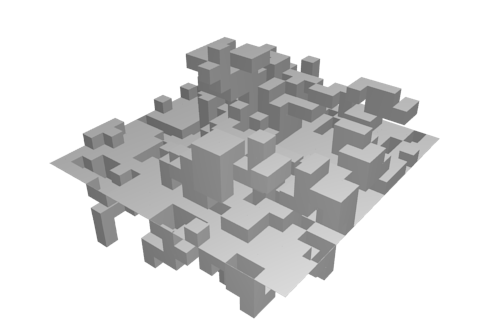} & \includegraphics[width = 0.5\textwidth]{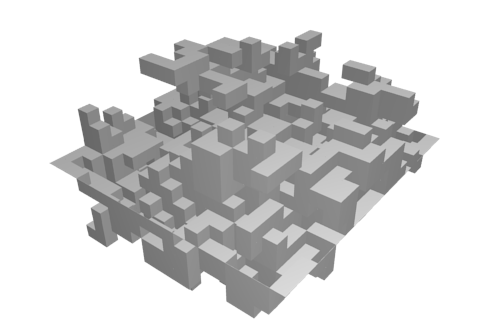} \\
\end{array}$
\caption{Several intermediate surfaces encountered when generating the LERS}
\end{figure}

For $n$ ranging from 5 to 100 we sampled from the uniform spanning 2-tree on $\mathcal{Q}_n$ 1000 times measuring the size of the loop-erased random surface with boundary the equatorial loop each time.

%%%%%%%%%%%%%%%%%%%%%%%%%%%%%%%%%%%%%%%%%%%%%%%%%%%%%%%%%%%%%%%%%%%%%%
\section{Results}
%%%%%%%%%%%%%%%%%%%%%%%%%%%%%%%%%%%%%%%%%%%%%%%%%%%%%%%%%%%%%%%%%%%%%%
%\subsection{present results and conclusions}
We have 1000 samples each for $n$ from 5 to 100. By our data, $\tilde{c} = 2.5269$.  In addition, we created 1000 bootstrap resamples of our data with a corresponding 1000 bootstrap estimates of $\tilde{c}$.  Discarding the lowest and highest 2.5\% gives us a nonparametric 95\% confidence interval for $c$ with lower bound 2.5252 and upper bound 2.5284.  That is, we are 95\% confident that $c = 2.5269 \pm 0.0017$. We expect a rational exponent similar to the $5/ 4$ for LERW in two dimensions and the conjectured $13/ 8$ for LERW in three dimensions \cite{Anton:2001}.  We hypothesize that, in fact, $c = 48/19$.  
\begin{figure}
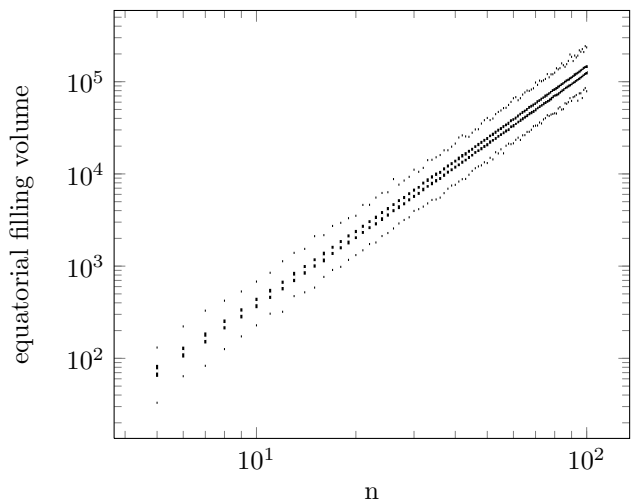

\mygraphone{}
\caption{For each $n$ the marks represent the typical data found in a bar and whisker plot}
\label{fig:eq}
\end{figure}

\section{Further Directions}
Collecting data for larger values of $n$ would help clarify the value of $c$ and either support or weaken the 48/19 hypothesis.  It is clear that the uniform spanning tree and loop-erased random walk in 3 dimensions are closely related to the loop-erased random surface.  What can be said about $c$ based on what is known and conjectured about these objects, specifically the conjectured growth exponent of $13/8$ for LERW in 3 dimensions.  Can it be proven that $c$ is strictly between 2 and 3?  That is, do there exist constants $\epsilon,c_1,c_2 >0$ such that as $n\to\infty$ we have $c_1 n^{2+\epsilon} < M_n < c_2 n^{3-\epsilon}$ with probability tending to 1? All these are interesting future directions for research into the loop-erased random surface in 3 dimensions.

\section{Acknowledgements}
I would like to thank Sayan Mukherjee and Randall Kamien whose comments on a draft of this paper were invaluable.  I would also like to thank my advisor Matthew Kahle, without whom none of this research would be possible.

\bibliographystyle{plain}
\bibliography{numexp_bib}

\end{document}